\documentclass{article}[12pt]
\usepackage{color,times,amsmath,amsfonts,latexsym,epsfig,epsf,colordvi}
\usepackage{url,hyperref}
\usepackage[english]{babel}

\newcommand{\CC}{\mathbb {C}}
\newcommand{\RR}{\mathbb {R}}

\newcommand{\bp}{\begin{pmat}}
\newcommand{\ep}{\end{pmat}}

\author{  Frank Uhlig \thanks{Department of Mathematics and Statistics, Auburn 
University, Auburn, AL 36849-5310 \ (uhligfd@auburn.edu)}}

 \title{\vspace*{-12mm} On the  Block-Decomposability of 1-Parameter Matrix Flows\\ and Static Matrices\\[-2mm]}
\begin{document}
\date{~}
\thispagestyle{empty}
\maketitle

\thispagestyle{empty}

\vspace*{-16mm}
\begin{center} { \bf Abstract  } \\[2mm]
\begin{minipage}{130mm}
For general complex or real 1-parameter matrix flows $A(t)_{n,n}$ and for time-invariant static matrices $A \in \CC_{n,n}$ alike, this  paper considers ways to decompose matrix flows and single matrices globally  via one constant matrix similarity $C_{n,n}$ as
$A(t)  =  C ^{-1} \cdot \text{ diag}(A_1(t), ..., A_\ell(t))   \cdot C$ or $A = C^{-1}\cdot  \text{diag}(A_1,...,A_\ell)\cdot  C$ 
with each diagonal block $A_k(t)$ or $A_k$ square and their number $\ell > 1$ if this is possible. The theory behind our proposed algorithm is elementary and uses the concept of invariant subspaces for the Matlab {\tt eig} computed 'eigenvectors' of one associated flow matrix $B(t_a)$ to find the coarsest simultaneous block structure for all flow matrices $B(t_b)$. The method works very efficiently  for all time-varying matrix flows, be they differentiable, continuous or discontinuous in $t$, and for all fixed entry matrices $A$;  as well as for all types of square matrix flows or fixed entry matrices such as hermitean, real symmetric, normal or general complex and real flows $A(t)$ or static matrices $A$,  with or without Jordan block structures and with or without repeated eigenvalues. Our intended aim is to discover diagonal-block decomposable flows as they originate in sensor driven outputs for time-varying matrix problems and thereby  help to reduce the complexities of their numerical treatments through adapting 'divide and conquer' methods for their diagonal sub-blocks. Our method is also applicable to standard fixed entry matrices of all structures and types. In the process we discover  and study k-normal fixed entry matrix classes that can be decomposed under unitary  similarities into various   $k$-dimensional block-diagonal forms.
\end{minipage}\\[-1mm]
\end{center}  
\thispagestyle{empty}


\noindent{\bf AMS :} 15A99, 15A21, 15B99,  65F99 
%\quad \ \ \ \ (ordered by significance)
\\[-7mm]

\pagestyle{myheadings}
\thispagestyle{plain}
\markboth{Frank Uhlig}{Decomposability of Matrix Flows and  Static Matrices }

\section*{Preface}\vspace*{-1mm} 
Matrix block decompositions have been studied for almost a century. They appeared first in the early days of quantum physics in the 1920s and were essential to comprehend how higher atomic weight elements might have been formed from lighter ones in the primeval stages of the universe under aggregation, pressure, and heat. These processes were modeled by 1-parameter hermitean matrix flows $A(t)$ and their  parameter-varying eigencurves were found to never cross for indecomposable hermitean matrix flows. This eventually  led matrix theoreticians and numericalists to study eigencurve crossings in order to grasp the notion of decomposable matrix flows. These studies have been conducted from the fixed entry, static matrix theoretical and  matrix numerics standpoints via backward stable methods and use matrix similarities, matrix factorizations, basic linear algebra subroutines (BLAs) and so forth; all in a Wilkinsonian way, paying attention to error analyses, to perturbation effects and so forth. Eigencurve crossings have been found and algorithms proposed to compute them, yet with little effect on actually decomposing decomposable hermitean or general complex matrix flows, nor how to tell whether a given flow is decomposable or not.\\[1mm]
This author has contributed to the eigencrossings literature in \cite{FUEigencross} but was ultimately not satisfied with his results. Soul searching, despair and small constructive steps then lead me to approach the matrix decomposability problem most simply via invariant eigenspace theory for flows and relying on logic 0-1 matrix structural computations. The resulting  elementary approach appears to be on the level of a master's thesis and not worth clogging up ten pages of any reputable Journal. \\[2mm]
But the method is new, it solves a previously almost intractable  matrix quandary and it does so very accurately and elegantly. It has many practical  applications. However, it  lies outside our current knowledge base for matrix theory and its applications. No error, no perturbation analyses are needed. And there are no AMS subject classifications for time-  or parameter-varying matrix flows, be they hermitean or general complex.  \\  
This is a new territory and  both, fixed entry matrices $A$ and general or hermitean matrix flows $A(t)$ can now be block-decomposed  very simply -- if that is  possible -- and new insights can be gained for static matrices and matrix flows alike. These can then be studied and computed in 'divide and conquer' fashion more speedily.\\[-6mm]

\section{Introduction }\vspace*{-1mm} 

This paper studies time-varying, i.e., 1-parameter varying matrix flows $A(t) \in \CC_{n,n}$ over an interval $t_o \leq t \leq t_f \in \RR$ or when $t$ follows a finite section of a curve in $\CC$ and fixed entry matrices $A \in \CC_{n,n}$. In many applications and in many matrix computations it is good to know whether a  dense matrix flow $A(t)$ or a dense general static matrix $A$  can be decomposed into an  array of diagonal blocks $A_k(t)$ or $A_k$ (for $k = 1,...,\ell$)  of smaller dimensions, i.e.,\\[1mm] 
\begin{equation} \label{blockdiag} A(t)_{n,n}  =  C ^{-1} \cdot  \bp 
A_1(t) & O & O &  \cdots & O\\
O& A_2(t)& O & \cdots & O\\   
O&O&\ddots &  \ddots & \vdots\\
\vdots & \vdots & \ddots & \ddots &O\\
O&O&\cdots & O & A_\ell(t) \ep  \cdot C \ \ \ \text{ or } \ \ \ A_{n,n} = C^{-1} \cdot \text{ diag}(A_1,...,A_\ell) \cdot C\ .
\end{equation}

\vspace*{1mm}
Here $C \in \CC_{n,n}$ is an  invertible fixed entry matrix that is invariant for all parameters $t$. If a matrix flow $A(t)$ or a static matrix $A$ can be  decomposed in this fashion,  then many numerical problems for $A(t)$ and $A$ may be 'divided and conquered' into $\ell$ smaller subproblems for the individual blocks $A_k(t)$ or $A_k$ and  these subproblems can usually be solved more quickly. The matrix flow $A(t)$  may derive over time from given equations or it may be generated from sensor data that arrives at a constant discrete sampling rate $\tau$ such as  $\tau = 0.02$ seconds  or $50 ~H\!z$.\\[2mm] 
%\newpage

Decomposable matrix flows have been intimately linked to eigencurve crossings of  matrix flows $A(t)$ for over 90 years. In 1927 and 1929, Hund \cite{FH1927} and  von Neumann and Wigner \cite{NW} proved  that hermitean single parameter matrix flows $A(t) = (A(t))^*$ whose eigencurves cross each other must be decomposable via a fixed unitary matrix $C$ in the above sense. An eigencurve crossing is sufficient for hermitean matrix flow decomposability, but  the converse implication is not true.  In \cite{FUEigencross}, the author has studied the eigencurves of hermitean and general matrix flows and developed an algorithm to deduct the coarsest block-diagonalization dimensions of hermitean  matrix flows from their eigencurve crossing data. The biggest drawback of that method for deciding matrix or data decomposability is the fact that  hermitean and general  matrix flows $A(t) \in \CC_{n,n}$ need not show eigencurve crossings at all, even if they are decomposable. This paper introduces a different algorithm that uses standard  matrix invariant subspace theory to decompose matrix flows into  block-diagonal flows -- if  possible  -- for both hermitean and general complex matrix flows, or it establishes that such decompositions are impossible for $A(t)$. These fundamental matrix flow results are then applied to the static matrix decomposability question and $k$-normal matrices in the Applications section.\\
Details, numerical codes and tests follow.\\[-7mm]

\section{Theory}\vspace*{-2mm}

This section deals with  elementary notions and facts for unitarily decomposable matrix flows $A(t) \in \CC_{n,n}$. Applications to static matrix  $A$ decompositions via unitary similarities are treated in the Applications section below.\\[1mm] 
To start we consider a 'proper' $n$ by $n$  hermitean time-varying matrix flow $A(t) \in \CC_{n,n}$ that can be block-diagonalized uniformly as described in (\ref{blockdiag}) for a fixed nonsingular matrix $C_{n,n}$  and all $t_o \leq t \leq t_f $. Let us assume that we have $\ell >1$ diagonal blocks here for 'properness'. Any hermitean matrix flow $A(t)$ allows us to diagonalize the flow matrix $A(t_a)$ for any $t_a \in [t_o, t_f]$ via a unitary similarity transformation $V(t_a)$ so that $A(t_a) \cdot V(t_a)  =  V(t_a) \cdot D(t_a)$ and $D(t_a)$ is real diagonal. The transforming  unitary matrix $V(t_a)$ contains the eigenvectors of $A(t_a)$ in its columns and the eigenvalues of $A(t_a)$ appear on the diagonal  of $D(t_a)$ if we use  Matlab's built-in  {\tt eig} function for example. As we have assumed that $A(t_a) = C^{-1} \cdot \text{blockdiag} (A_1(t_a), ..., A_\ell(t_a)) \cdot C$ for  some nonsingular fixed entry matrix $C$, each  eigenvector in $V(t_a)$ is associated with one of the eigenvalues of $A(t_a)$ and in fact the eigenvector columns of $V(t_a)$ that are associated with the eigenvalues of one diagonal block $A_i(t_a)$ form an orthonormal  basis for an invariant subspace of  $A(t)$ of which there are $\ell$ by assumption. \\[1mm]
As all matrices of our assumed decomposable  hermitean flow $A(t)$ share the same invariant subspace structure expressed in (\ref{blockdiag}), then for any $t_b \neq t_a \in [t_o, t_f]$ the matrix $A(t_b)$ must be block diagonalizable and \\[-3mm]
\begin{equation} \label{two}
\tilde A_{t_a}(t_b) = (\tilde V(t_a))^* \cdot A(t_b)   \cdot \tilde V(t_a)
\end{equation}

\vspace*{-1mm}
 will be block-diagonal with the  same common block structure as soon as  we have re-arranged the eigenvector columns of $V(t_a)$ in $\tilde V(t_a)$ into $\ell$  groups that generate equal zero and non-zero pattern columns in $ \tilde A_{t_a}(t_b)$. The re-arrangement of the columns of $V(t_a)$ can be achieved by looking at the logic 0-1 {\tt spy} matrix $Al(t_b)$ of $\tilde A_{t_a}(t_b)$ in Matlab when all entries in $\tilde A_{t_a}(t_b)$ below a certain magnitude threshold have been set equal to zero.  Then we re-sort the columns of $V(t_a)$ so that  logic 0-1 vectors in $Al(t_b)$ (and thus of $\tilde A_{t_a}(t_b)$) fall  into $\ell$ distinct groups according to the location of their almost-zero and non-zero entries and thereby obtain $\tilde V = \tilde V(t_a)$.\\[1mm]
{\bf Theorem 1}: \begin{minipage}[t]{142mm}{If a hermitean time-varying matrix flow $A(t)$ can be properly and  uniformly diagonalized by a constant unitary similarity $U^* ... \ U$, then the eigenvector matrix $V$  of any flow matrix $A(t_a)$  can be re-arranged column-wise in $\tilde V$ so that  any  matrix $\tilde V^*A(t_b)  \tilde V$ with $t_b \neq t_a$ has the identical or a finer block-diagonal structure for all $t_{b}$. And vice versa, if the eigenvector matrix $V$ of one hermitean matrix flow matrix $A(t_a)$ creates a block-diagonalizable logic 0-1  pattern matrix for $A(t_b)$ with $t_b \neq t_a$ under a column-rearranged version $\tilde V$ of $V$, then all matrices in the hermitean flow $A(t)$ are simultaneously block-diagonalizable by the same unitary similarity  $\tilde V^* \ A(..) \ \tilde V$.}\\[-2mm]
\end{minipage}
For general complex matrix flows $A(t)_{n,n}$ that are diagonalizable throughout their single-parameter range, the same invariant subspace argument holds except that the unitary eigenvector matrix $\tilde V(t_a)$ similarity needs to be replaced by a general similarity via a nonsingular matrix $\tilde W(t_a)$ so that the inverse $\tilde V(t_a)^*$ of $\tilde V(t_a)$  in formula (\ref{two}) becomes $\tilde W(t_a)^{-1}$.\\[1mm]
{\bf Theorem 2}: \begin{minipage}[t]{142mm}{If a diagonalizable general complex time-varying matrix flow $A(t)$ can be properly and  uniformly diagonalized by a constant matrix similarity $C^{-1} ...\  C$, then the eigenvectors of any flow matrix $A(t_a)$ will block-diagonalize -- upon re-sorting -- every other flow matrix $A(t_b)$ by similarity into block-diagonal  form which may be finer than the coarsest possible block-diagonal form of that flow. And vice versa.}\\[-2mm]
\end{minipage}
Here the term 'coarsest block-diagonal form' refers to one with the minimal possible block number $\ell$ in (\ref{blockdiag}). Note for example, that $D(t_a)$ in formula (\ref{two})  represents the finest, i.e., a 1 by 1 block-diagonalization with  $\ell = n$   for $A(t_a)$.\\[2mm]
The next section deals with computing the minimal number $\ell$ of invariant subspaces of a properly decomposable matrix flow $A(t)$ by re-sorting the columns of their respective eigenvector matrix $V(t_a)$ or $W(t_a)$, so that the coarsest  simultaneous diagonal block reduction (or a finer one) can be achieved for any flow matrix $A(t)$ effectively -- provided that $\ell$ is found to exceed 1.\\[-9mm] 

\section{The Algorithm and Computed Results}\vspace*{-2mm}

As   theory tells us, to solve the matrix flow decomposability problem it suffices to compute\\[1mm]
 ({\bf A}) the eigenvector matrix $X(t_a)$  of one  flow matrix $A(t_a)$ and apply the  similarity $(X(t_a))^{-1} \cdot  \ A(t_b) \ \cdot X(t_a)$ to\\
 \hspace*{5mm} any other flow matrix $A(t_b)$\\[1mm]
  in order to learn about the coarsest (or a  finer) block-diagonalization of the given matrix flow. Theory predicts perfect zeros in the updated $(X..)^{-1} \cdot  A..  \cdot X..$ flow matrix, but numerical rounding errors and conditioning problems  always create relatively small magnitude entries in some entry positions of  the computed $\hat A(t_b) = (X(t_a))^{-1} \cdot A(t_b) \cdot X(t_a)$ that theoretically ought to be zero. These tiny magnitude entries must be replaced by zeros in order to exibit the block structure  of $A(t_b)$ properly.\\
  For this purpose we form \\[1mm]
 ({\bf B})  the logical 0-1 matrix of the computed $\hat A(t_b)$ matrix  and then  \\[1mm]
 ({\bf C}) we rearrange its rows and columns by collecting equal 0-1 pattern row vectors therein into groups in order to\\
 \hspace*{5mm} exhibit the block-diagonal structure of the studied flow. \\[1mm]
This process works  equally well for all time-varying matrix flows. It offers a great improvement over what could be gleamed geometrically from eigencurve crossings in \cite{FUEigencross}. Besides, in  \cite{FUEigencross} the general complex matrix flow case was generally  found to be intractable via  coalescing eigencurve studies. Here this problem does not even appear.\\[2mm]
Figure 1 shows nine Matlab {\tt spy} graphs for a dense non-normal complex 17 by 17 matrix flow $A(t)$. Reading  this figure  row by row, the first row of graphs shows the Matlab  {\tt spy} 0-1 pattern transitions from $A(t_a)$ to $A(t_b)$; the second row shows the  0-1 pattern transitions from $A(t_b)$ to $A(t_c)$, and the third row  shows the ones from $A(t_c)$ to  $A(t_{rd})$ for a randomly chosen parameter $t_{rd} \in \CC$. Column {\bf (A)} shows the diagonalization $D(t_a)$ 0-1 pattern for $A(t_a)$ via Matlab's {\tt eig} function. Column {\bf (B)} displays the  0-1 similarity patterns from varying starting matrices $A(t_a), ..., A(t_c)$. Note that  the spy graphs in column {\bf (B)} all hint at a 7 4 3 2 1 block-diagonalization for this flow. The third column {\bf (C)} of  {\tt spy} graphs  is computed from the  0-1 data in column {\bf (B)} by collecting equal 0-1 row vectors in groups.  In {\bf (C)}  the same non-zero diagonal blocks of dimensions 7,  4, 3, 2, and 1 appear, but they are arranged in differing orders.\\[1mm]
 The general complex matrix flow $A(t)_{17,17}$ of this test example was built from a  complex matrix flow  $B(t)$ with   block-diagonal dimensions 7 4 3 2 1. $B(t)$  was then  transformed into the dense general flow $A(t) \in \CC_{17,17}$ by a fixed random entry dense unitary similarity. The spy graph sequences in Figure 1 below were  computed by our MATLAB algorithm {\tt deccomplflow9.m} in the subfolder \emph{general flows} of \cite{FUMatrixDecomp}.\\[1mm]
% \newpage
 \hspace*{21mm} ({\bf A}) \hspace*{46mm} ({\bf B}) \hspace*{45mm} ({\bf C}) \\[-6mm]
\begin{center}
\includegraphics[width=152mm]{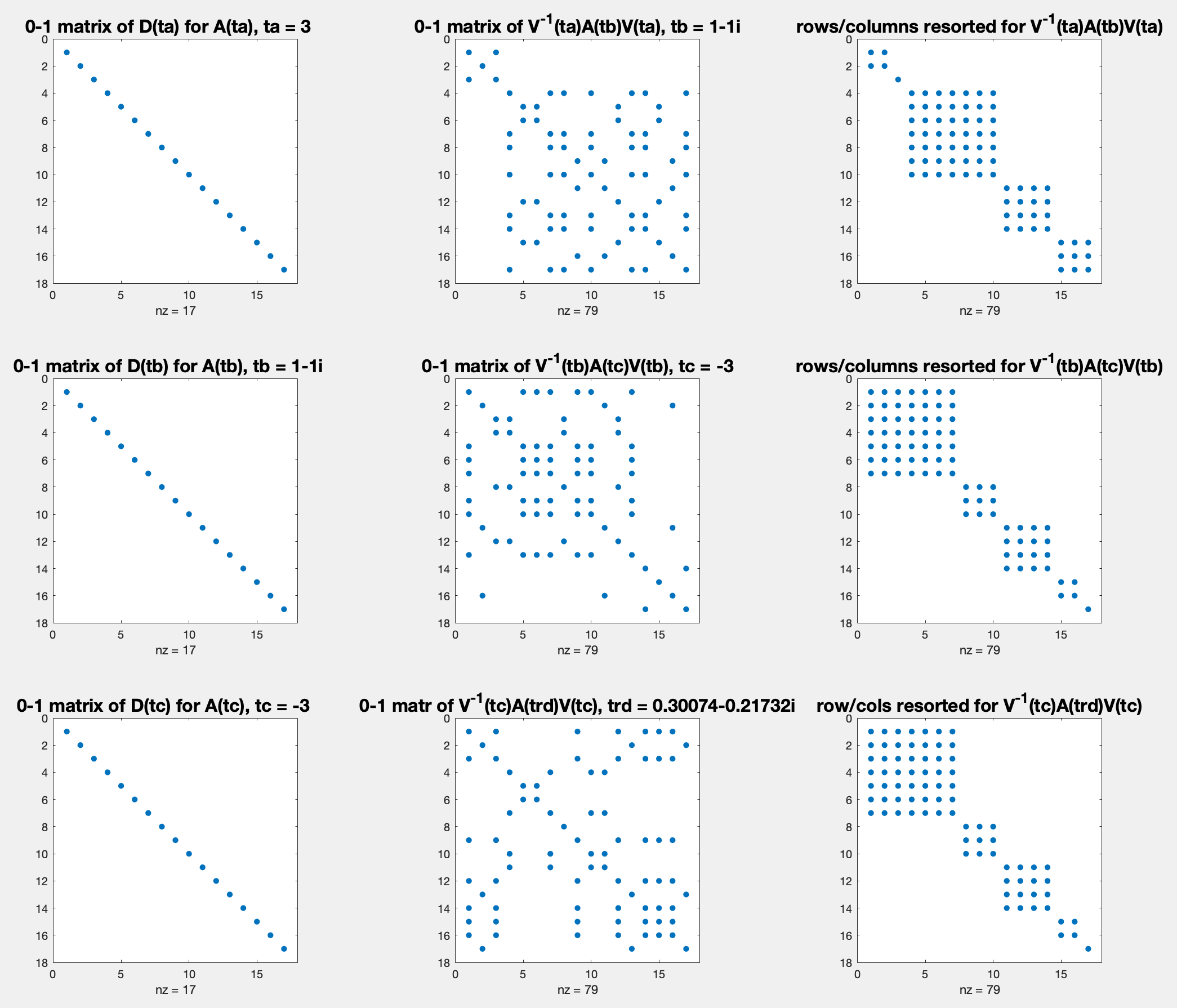} \\[0mm]
Figure  1: General complex flow $A(t)$ with a chain progression from $t = t_a$ to $t_b$ to $t_c$, and to random $t=t_{rd}\in \CC$
\end{center}

Our algorithm works equally well for matrix flows that are built from proper Jordan blocks such  as the next 9 by 9 complex  flow  example shows  with Jordan blocks of sizes 4 and 5 in  Figure 2. Note the 'holes with zeros' in the respective diagonal 0-1 {\tt spy} blocks in columns {\bf (B)} and {\bf (C)} that seem to occur occasionally for Jordan blocks. \\
We know not why.\\[1mm]
%\newpage
\hspace*{23mm} ({\bf A}) \hspace*{46mm} ({\bf B}) \hspace*{46mm} ({\bf C}) \\[-6mm]
\begin{center}
\includegraphics[width=155mm]{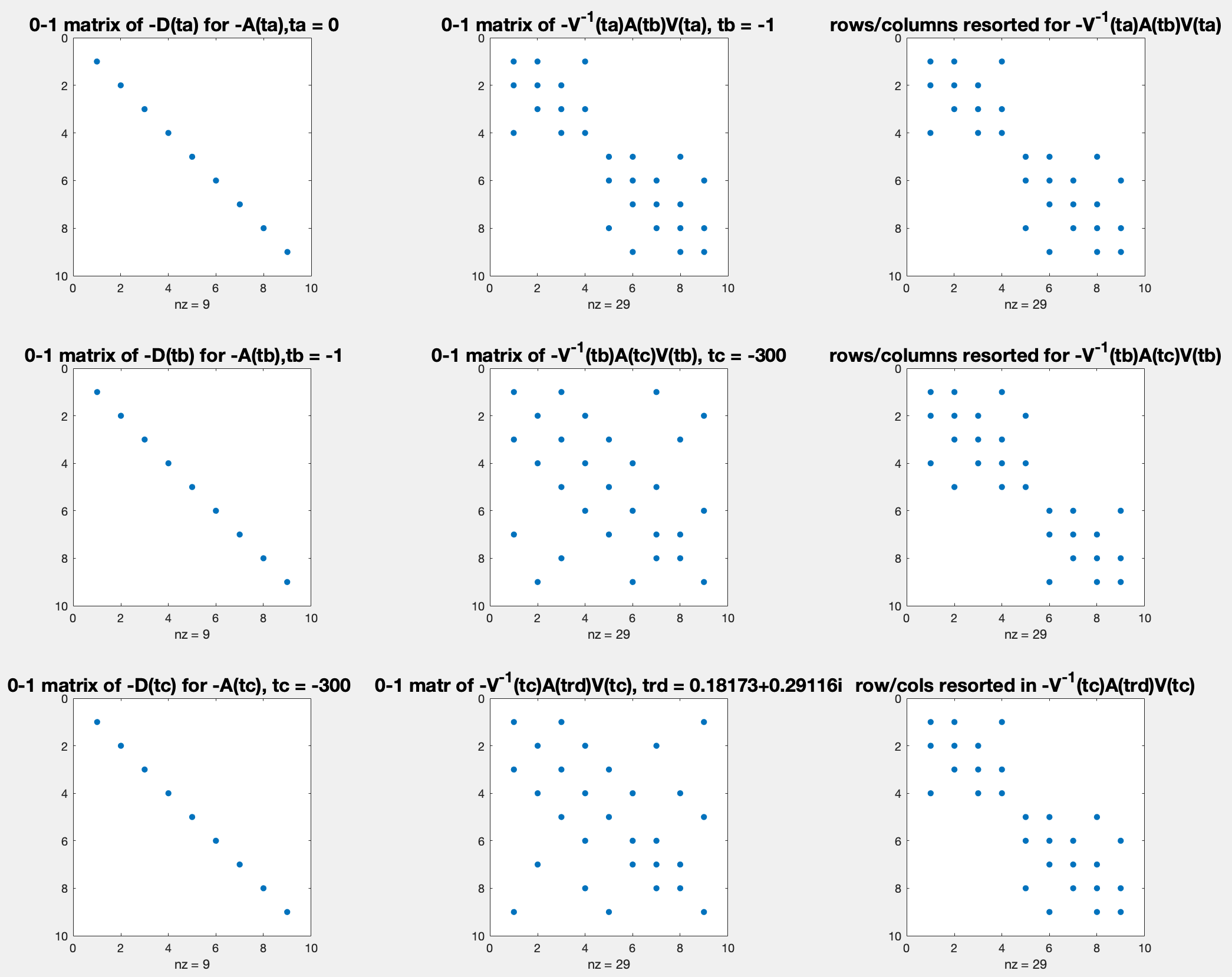} \\[0mm]
Figure  2: Dense decomposable complex matrix flow $A(t)$, formed from two Jordan blocks of sizes 4 and 5
\end{center}
By depending  only on elementary invariant subspace theory, our algorithm and code works well with real time parameters $t \in \RR$, as well as with more general complex parameters $t \in \CC$ as shown in  Figures 1 through 4.\\[1mm]
The web depository \cite{FUMatrixDecomp} also contains a simpler algorithm ({\tt deccompl.m}) for  finding the block-diagonal dimensions of a general matrix flow for just one time $t_b$ from the Matlab {\tt eig} diagonalization of $A(t_a)$ with $t_a \neq t_b$. Besides, there is a different 9-graph Matlab m-file ({\tt deccomplflow9a.m}) in \cite{FUMatrixDecomp}  that computes the pattern transitions not along the chain from $t_a \text{ to } t_b, \text{ then to } t_c$ and then from $ t_c \text{ to } t_{rd}$ as {\tt deccomplflow9.m} does, but instead  computes the transitions starting always from $t_a$ to each of $t_b, t_c \text{ and } t_{rd}$ in turn.\\
 The \emph{Matrixflow Decomp} folder at \cite{FUMatrixDecomp}  contains the  Matlab m-files  for general 1-parameter matrix flows in  the subfolder \emph{general flows}. The subfolder \emph{hermitean flows}  at \cite{FUMatrixDecomp} deals with hermitean or symmetric single-parameter matrix flows. The hermitean flow methods  {\tt  decherm.m, dechermflow9.m} and {\tt dechermflow9a.m} are made simpler by the fact that they do not have to deal with the Matlab {\tt eig.m} m-file output for derogatory non-normal matrices with proper Jordan block structures. There is also a subfolder \emph{staticMatrdecomp} in \cite{FUMatrixDecomp} with a block decomposition code for decomposable static matrices $A \in \CC_{n,n}$. This is used and explained in the Applications section.\\[2mm] 
The occurrence of Jordan blocks in a general complex matrix flow $A(t)$ and their treatment in {\tt eig.m} may also create bands of  0-1 entry rows of all 1s when computing $\hat A(t_b) = (X(t_a))^{-1}  A(t_b)  X(t_a)$. The all-1s rows need to be taken care of differently in the general case than in the hermitean matrix flow case, where such can not occur.\\[1mm]
In column {\bf (B)}, Figure 3 below shows such a banded  0-1 pattern matrix with several all-1s rows for a dense  example flow $A(t)_{14,14}$ that was built from a general complex 14 by 14 matrix flow that contains two Jordan blocks of size 2 and other blocks of dimensions 1 (3-fold), 3, and 4. Note that grouping identical 0-1 rows of column ({\bf B}) together into one diagonal block in column ({\bf C}) for Jordan block containing general flows -- as is sufficient for hermitean flows -- would result in all  rows of the 0-1 pattern matrix  becoming  indistinguishable here, indicating falsely that this general complex flow is indecomposable. \\
 For general non-hermitean flows  the actual re-sorting from ({\bf B}) to ({\bf C}) {\tt spy}  matrices uses both the zero and the non-zero pattern of each not-all-1s row of a {\tt spy} graph in column ({\bf B}) to arrive at the 0-1 {\tt spy} graph in column ({\bf C}). This helps us detect the block-diagonal dimensions correctly while also allowing us to determine  the total sum of all Jordan block dimensions for such flows.\\[1mm]
Each of our Matlab codes  provides on-screen interpretations of the computed outputs and describes the resulting block dimension sizes, for both hermitean and general complex matrix flows. For the latter, the on-screen block dimensions  refer to the summed dimensions of all Jordan blocks in the listed flow dimensions if followed by a (J). On-screen, there are warnings when the norm of an intermediate matrix $A(t...)$ becomes excessively large in which case the computed block dimension results may be erroneous or unreliable.\\[1mm]
\hspace*{19mm} ({\bf A}) \hspace*{48.5mm} ({\bf B}) \hspace*{49mm} ({\bf C}) \\[-6mm]
\begin{center}
\includegraphics[width=155mm]{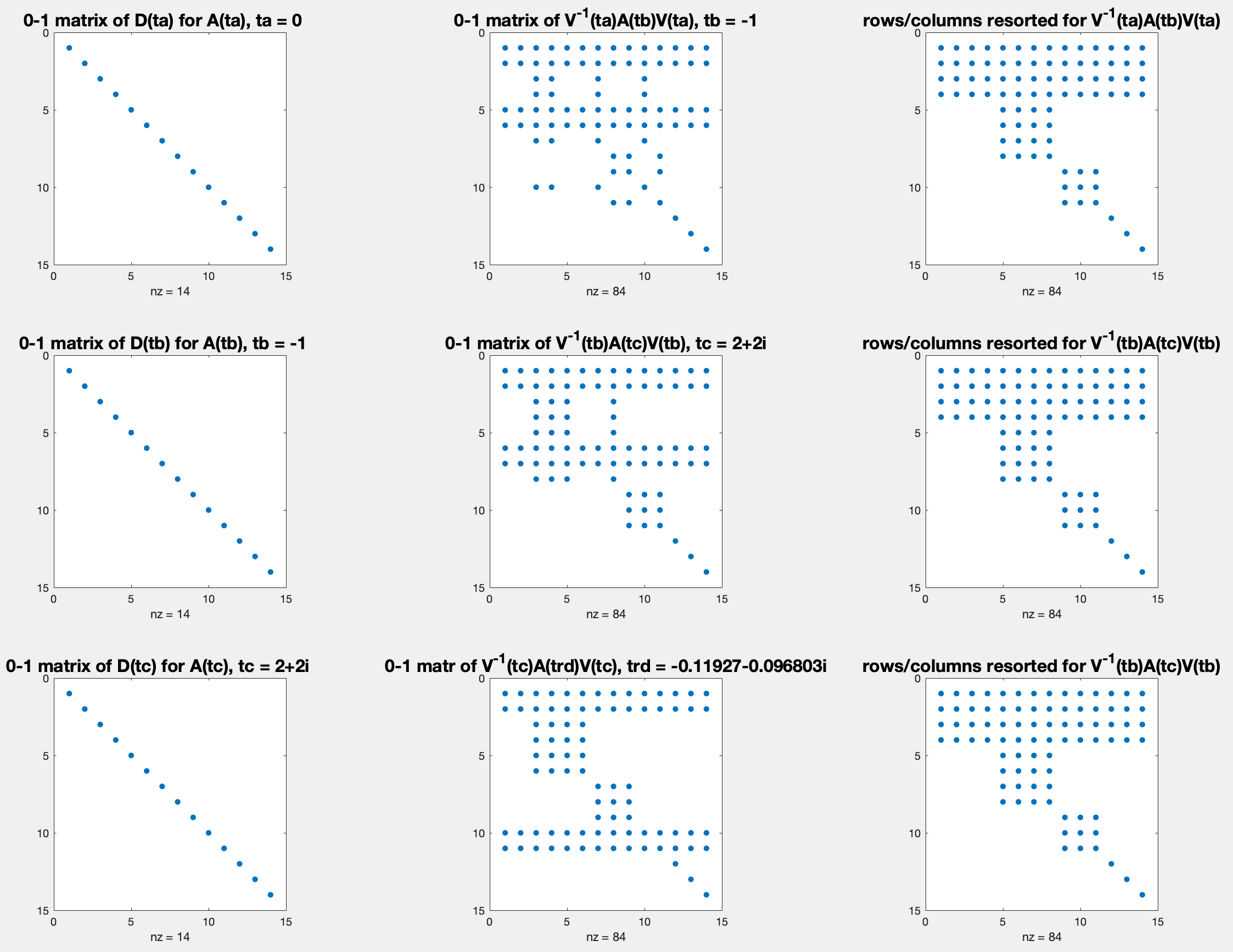} \\[0mm]
Figure  3: General decomposable complex matrix flow $A(t)$ with two 2 by 2 Jordan blocks
\end{center}
Note that in the intermediate second column the graphical 0-1 output ({\bf B})  in  Figures 1 through 3 may  differ from row to row, but that the final graphs in column ({\bf C}) are identical in each test run, except for possible permutations in the  order of the diagonal blocks.\\[1mm]
The use and success of our matrix flow decomposition algorithm does not depend on or require any smoothness conditions for any given  1-parameter matrix flow $A(t)$. Real-time or discontinuous data from sensor inputs is quite admissible. Figure 4 below shows the block pattern output for a general complex matrix flow $A(t)$, constructed  from a block diagonal time-varying matrix flow $B(t)$ whose blocks have uniform block sizes for all parameters $t$. In $B(t)$ some diagonal blocks have  mixed time-varying entries and random entries that change erratically with every call. The matrix flow $B(t)$ is then made dense  to become $A(t) =  V^*B(t)V$ by using  a  randomized but fixed unitary matrix similarity $V$  to create  a dense test example with partial random entries that is  decomposable.\\[1mm]
\hspace*{19mm} ({\bf A}) \hspace*{48.5mm} ({\bf B}) \hspace*{49mm} ({\bf C}) \\[-6mm]
\begin{center}
\includegraphics[width=155mm]{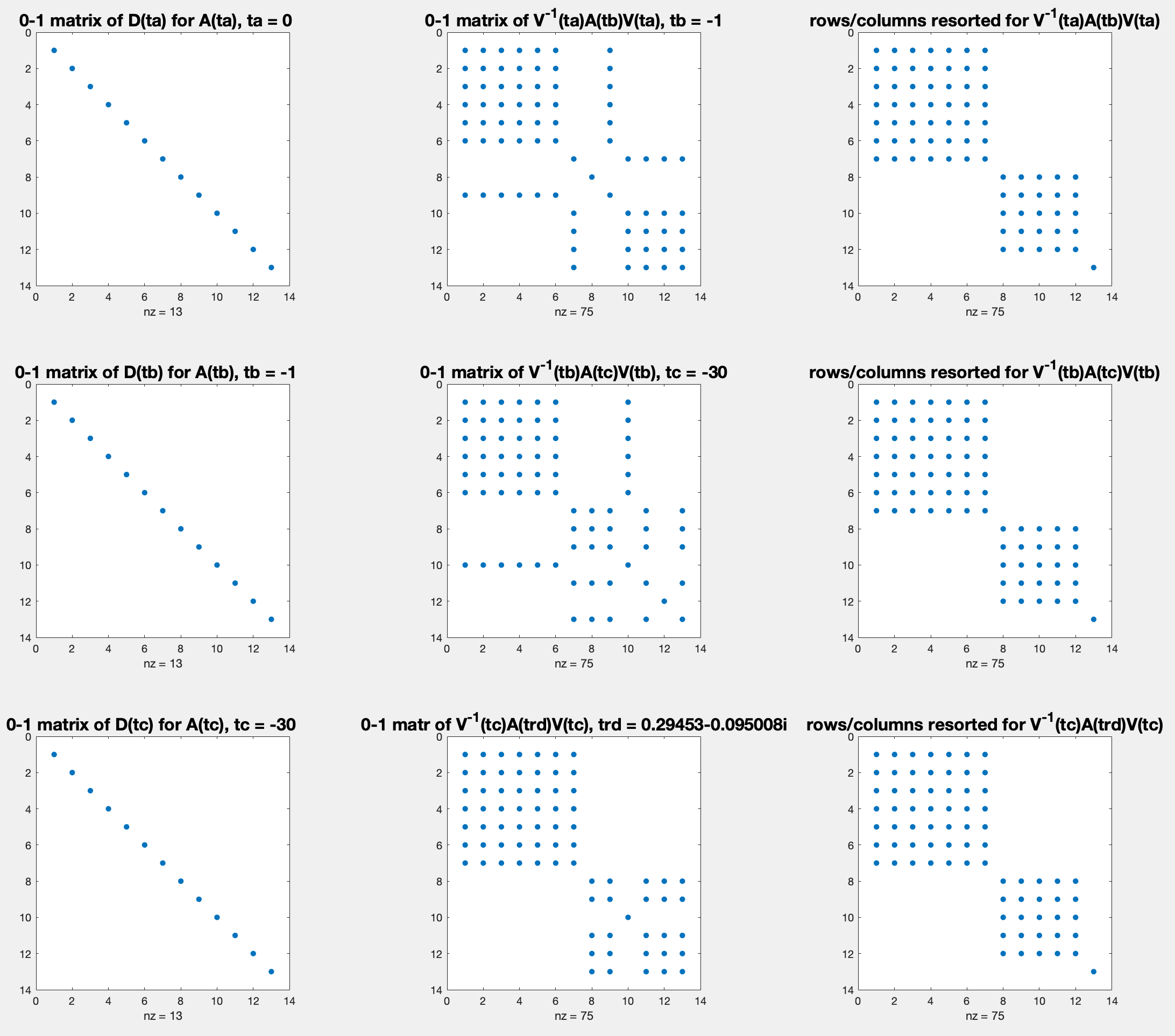} \\[0mm]
Figure  4: Decomposable complex matrix flow $A(t)$ with  7-, 5- and 1-dimensional  random entry diagonal blocks
\end{center}

Our matrix flow decomposition algorithms take very little time on a 2019 MacBook Pro, around 0.04 sec for matrix flows $A(t)_{n,n}$  of size 250 by 250 and about 0.9 sec when $n = 1000$. Overall they require one eigenanalysis of $A(t_a)$, several matrix multiplications and some  logic matrix arithmetic. \\[1mm]
The depository \cite{FUMatrixDecomp} includes  Matlab codes for constructing more than a dozen example flows in both, the hermitean and the general case. Test problem generation  can be implemented inside our respective flow decomposing routines by entering an integer matrix dimension number $n \leq 20$. Plotting can be turned off by setting the input parameter {\tt zeich} unequal to 1 and the block dimension information will still be displayed on-screen. Without graphing, the CPU times for running the algorithms were in the hundreds of a second in all tested dimensions and -- with graphing included -- the computations and all visual '{\tt spy}' displays would appear in a fraction of a second.%\\[-6mm]

\section {Applications}\vspace*{-2mm}

When used in a 'preconditioner' type of way, any matrix block decomposition algorithm (${\cal MBD\!A}$) may be of benefit in many matrix computational problems, both for fixed entry matrices $A \in \CC_{n,n}$ and for  time-varying matrix flows $A(t) \in \CC_{n,n}$. In particular if an ${\cal MBD\!A}$ is applied to a similarity invariant or unitarily invariant matrix  problem and the given matrix $A$ or matrix flow $A(t)$ are decomposable, then  the CPU time savings for computing the desired properties of $A$ or $A(t) \in \CC_{n,n} $ amount to around 50 \% of CPU time if the largest diagonal block has size $m = 0.8 \cdot n$ after ${\cal MBD\!A}$.  And if the largest diagonal block dimension after a ${\cal MBD\!A}$ reduction has  size $m = n/2$, then the  savings would reduce the original CPU time by at least  three quarters. Here we assume that an $O(n^3)$ process is being used on the flow's diagonal blocks in a 'divide and conquer' algorithm after an ${\cal MBD\!A}$ preconditioning.  These savings rates have been established in \cite{FUDecompMatrFoV}. In \cite{FUDecompMatrFoV} the field of values problem for a fixed entry matrix $A$ was studied in light of  $A$'s decomposability and our elementary block diagonalization method was applied to the hermitean matrix flow $H\!K_A(t) = \cos(t) H + \sin(t)K$ with $H = (A+A^*)/2$ and $K = (A-A^*)/(2i)$. This  achieved   speed-ups of up to 12 times when compared with the standard QR based matrix eigen-algorithm  for computing the field of values boundary curve accurately for decomposable matrices $A$; again see \cite{FUDecompMatrFoV}.\\[1mm]
Note also that the  field of values idea can  be used in reverse to find the eigenvalues of decomposable general fixed entry matrices $A \in \CC_{n,n}$ with less effort. To find a possible block diagonalization of any square matrix $A\in \CC_{n,n}$ by our elementary ${\cal MBD\!A}$ method requires us to compute the eigenvalues of one hermitean matrix $H\!K_A(t_a)$ that is derived from $H = (A+A^*)/2$ and $K = (A-A^*)/(2i)$ for $A$. This can be done  at much smaller hermitean $O(n^3)$ QR cost than a general matrix QR eigenanalysis on $A$ itself. Thereafter we check the 0-1 pattern matrix of any $H\!K_A(t_b)$ where $t_b \neq t_a$. If its logic 0-1 pattern shows a  block diagonalisable pattern we re-sort the  eigenvector columns for $H\!K(t_a)$ that reside in $V$ into  equal row pattern groups in $\tilde V$ and then the matrix $\tilde A = (\tilde V)^* \cdot A \cdot \tilde V$ will be block diagonal per Theorems 1 and 2.  $\tilde A$'s eigenvalues can then be found more expediently from the smaller general diagonal blocks of $\tilde A$ than using Matlab's {\tt eig} function on the originally dense $n$ by $n$ matrix $A$.\\
 The overall cost of finding the eigenvalues for  a general dense, but decomposable matrix $A$ thus is essentially reduced to the cost of one hermitean $n$ by $n$ matrix eigenanalysis plus several smaller sized general matrix block eigenanalyses and very cheap logic matrix overhead.\\[2mm]
 Here is one previously unknown  matrix theoretical application: We have tested several,  mostly  non-normal 'gallery' test matrices of Matlab for block diagonalisability and found three non-normal matrix classes that are unitarily similar to 2 by 2 or 4 by 4 block diagonal matrices, as well as one dense static $n$ by $n$ matrix type that can be decomposed  into exactly two almost equal sized diagonal blocks for all $n$. The  concept of block diagonalisability introduces a new concept here. It allows us to generalize standard normal matrices $A$ with $A^*A = AA^*$. Normal matrices  can rightfully be called 1-normal matrices  since they  can  always be unitarily diagonalized into 1 by 1 block diagonal form. Next we define define $k$-normal matrices.\\[1mm]
 {\bf  Definition} : \emph{A matrix $A \in \CC_{n,n}$ is called {\bf $\mathbf k$-normal} if A can be unitarily diagonalized  into block diagonal form with the maximal number of blocks of size $k$ that fits into $n$ and  smaller  dimensioned blocks depending upon the divisibility of $n$ by $k$.}\\[1mm]
The Matlab 'gallery' matrices {\tt 'binomial'}, {\tt "clement'}, and {\tt 'invol'} are non-normal. The {\tt 'invol'}  test  matrices with dimensions $2 \leq n \leq 14$ can be verified  to be 2-normal by using our ${\cal MBD\!A}$. 
The matrix norms of higher dimensional {\tt 'invol'} matrices reach astronomical heights; for example for $ n= 200$ the norm of {\tt 'invol(200)'} has order $10^{212}$ with essential entries of magnitudes between $10^{200+}$ and 1. For $n = 15$ the {\tt 'invol'} matrix norm exceeds $10^{11}$ and our algorithm can no longer differentiate between its huge and  rather small, but significant entries when we use fixed  rounding error thresholds when forming  logic 0-1 {\tt spy} matrices.\\[0.6mm]
 Matlab's {\tt 'binomial'} matrices are chameleon-like here: if $n$ is divisible by 4, {\tt 'binomial'} matrices decompose under unitary similarity into $n/4$ four  by four diagonal  blocks. For even dimensions $n$ that are not divisible by 4, {\tt 'binomial'} matrices diagonalize into all 4 by 4 diagonal blocks with  two additional 1 by 1 diagonal blocks. For odd $n = 2j+1$, however,  all {\tt 'binomial'} $n$ by $n$ matrices are 2-normal with $j$ two by two diagonal blocks and one additional 1-dimensional diagonal block. These results were obtained for $3 \leq n \leq 30$. At $n = 31$, the norm of 'binomial(31)' exceeds $10^8$ and our fixed threshold computations stop making sense.\\[0.6mm]
Quite differently again, the non-normal Matlab gallery matrix  {\tt 'clement'}  can be unitarily block reduced to two almost equal sized diagonal blocks for any dimension $n$.  This  makes {\tt 'clement'} $n$ by $n$ matrices  ${[}(n+1)/2{]}$- normal. Here the symbol ${[} .. {]}$ denotes the greatest integer function.
 The static {\tt 'clement'} matrix block decompositions have been achieved successfully by our ${\cal MBD\!A}$ code for all $n \leq 200$. Note that the norms of  {\tt 'clement(n)')} matrices stay uniformly well below 10 for $n \leq 200$ and ${\cal MBD\!A}$ computations for {\tt 'clement'} matrices were not tried for larger than dimension $n =  200$.\\[0.6mm]
Finally we noted that the non-hermitean {\tt 'circul'} matrices of Matlab were diagonalized by our ${\cal MBD\!A}$, making the {\tt 'circul} normal which we then checked via the $A^*A = AA^*$ equation, but did not realize before.\\[1mm]
 We do not know if these specific matrix diagonal block reducibilities via unitary similarities  are known or not, nor whether they have ever been exploited. We know how to exploit standard 1-normality for unitarily invariant static matrix problems in their computations; why not try to develop efficient methods for 2-normal, 3-normal or 4-normal static matrices as well.\\[1mm]
 Our ${\cal MBD\!A}$ Matlab code {\tt decompstaticMatr.m} is in the subfolder {\em staticMatrdecomp} of \cite{FUMatrixDecomp}. The code requires one input, a general static matrix $A \in \CC_{n,n}$. It has three outputs, a unitarily similar block diagonal matrix $Ad$ of $A$ if $A$ is decomposable, the unitary transforming matrix $V\!c$, and a list of the  block dimensions of $Ad$ if $A$ was decomposed. On-screen comments explain the results.  A run with  Matlab gallery matrices that were exemplified above and small $ n \leq 30$ (when feasible)  takes around 0.25 seconds of CPU time and calls with $n = 200$ when feasible take around 0.34 seconds on a 2019 MacBook Pro.\\[-5mm]

\section{An Outlook and Adjacent  Areas of Research}\vspace*{-0mm}

It might be of interest to size each occurring Jordan block in a general matrix flow $A(t)$ individually in our Figures 3 and 4 in the column {\bf B} rather than summarily, but we have not done so. Regarding  Jordan structures of fixed  static matrices $A_{n,n}$, it appears to be nearly impossible in general and at least very expensive to try and determine the Jordan structure of even small dimensioned static matrices $A$ reliably by numerical means such as {\tt eig} in Matlab. More involved computational efforts to find the Jordan normal form of small static matrices $A_{n,n}$ reliably are  in \cite{TE05} and similar efforts for the Kronecker normal form of singular matrix pencils are in \cite{EK17}. Yet the problem of block-diagonalizing  time-varying general matrix flows $A(t)$ in the presence of  Jordan structures  or  of repeated eigenvalues has been easily answered  computationally  here by using elementary invariant subspace theory. This shows that time- or single parameter-varying matrix flows $A(t)$  follow different fundamental concepts than classic static matrix theory and analysis. \\[1mm]
Could one and how could one alter the Francis multishift implicit QR method, for example, to account for repeated eigenvalues and higher dimensional principal  subspaces of static matrices somehow, we wonder.\\[1mm]
An application of our matrix flow decompositions helps with the matrix field of values problem for decomposing  general static  complex  matrices $A_{n,n}$, see  \cite{FUDecompMatrFoV}. Our matrix flow decomposition algorithm   now allows path following methods to compute the field of values of such matrices more efficiently than global matrix eigensolvers such as QR. \\[2mm]
Separately Loisel and Maxwell \cite[Thm 2.5, Sect 5, 6.2, and 7.1]{LM18} have developed an IVP ODE solver to find eigencrossing points of hermitean block-diagonal matrix flows for the field of values boundary computation problem, while Dieci et al \cite{DE99, DPP13,DP08,DGPP11, DP12} have studied multi-parameter flows and their eigencrossings as well as singular value crossings using geometric localization and zoom-in optimization methods. Maybe our invariant subspace based idea can be extended and adapted to help with such problems.\\[1mm]
Finally, Sabuya \cite{S65} has dealt with a related problem to classify all matrix flows $A(t)$ that are block-diagonalizable under time-varying similarities $X^{-1}(t) \cdot A(t) \cdot X(t)$ in contrast to our unified fixed $C ^{-1}\cdot A(t) \cdot C$ block-diagonal similarities.\\[-5mm]

\section*{Acknowledgement}
\vspace*{3mm}

\noindent
\centerline{{[} .. /box/local/latex/Decompmatrixflow.tex] \quad \today }

\vspace*{2mm}

\noindent
 4 image files :\\[2mm]
Decgen9pics17.png\\
Decgen9apics9.png\\
Decgen9pics14.png\\
Decgen9pics13.png\\

\end{document}